\documentclass[12pt]{article}
\usepackage{amsmath,amssymb,amsthm,amscd,graphicx}
\usepackage{latexsym}
\usepackage{enumerate}
\usepackage{amsbsy, amsfonts, textcomp}
\usepackage{makeidx}
\newtheorem{theorem}{Theorem}[]
\newtheorem{proposition}[theorem]{Proposition}

\newtheorem{corollary}[theorem]{Corollary}

\theoremstyle{definition}
\newtheorem{definition}[theorem]{Definition}
\newtheorem{example}[theorem]{Example}

\theoremstyle{remark}
\newtheorem{remark}[theorem]{Remark}
\newcommand{\F}{\mathcal{F}}
\newcommand{\K}{\mathcal{K}}
\newcommand{\G}{\mathcal{G}}
\newcommand{\M}{\mathfrak{M}}
\newcommand{\cL}{\mathcal{L}}
\newcommand{\R}{\mathbb{R}}
\newcommand{\Q}{\mathbb{Q}}
\newcommand{\C}{\mathbb{C}}
\def \GL{\operatorname{GL}}
\def \Gal{\operatorname{Gal}}
\def \SO{\operatorname{SO}}
\def \Ker{\operatorname{Ker}}
\def \DGal{\operatorname{DGal}}

\begin{document}

\Large
\centerline{\bf Real Picard-Vessiot theory}

\large
\centerline{Teresa Crespo, Zbigniew Hajto, El\.zbieta Sowa}
\normalsize
\let\thefootnote\relax\footnotetext{T. Crespo and Z. Hajto acknowledge support of grant MTM2009-07024, Spanish Science Ministry.}
\begin{abstract}
The existence of a Picard-Vessiot extension for a homogeneous linear differential equation has been established when the differential field over which the equation is defined has an algebraically closed field of constants. In this paper, we prove the existence of a Picard-Vessiot extension for a homogeneous linear differential equation defined over a real differential field $K$ with real closed field of constants. We give an adequate definition of the differential Galois group of a Picard-Vessiot extension of a real differential field with real closed field of constants and we prove a Galois correspondence theorem for such a Picard-Vessiot extension.
\end{abstract}

\tableofcontents

\section{Introduction}

Picard-Vessiot theory denotes Galois theory of homogeneous linear differential equations. The Picard-Vessiot extension associated to a given homogeneous linear differential operator is the analog of the splitting field of a given polynomial. For a homogeneous linear differential equation $\mathcal{L}(Y)=0$ defined over a differential field $K$ with field of constants $C$, a Picard-Vessiot extension is a differential field $L$, differentially generated over $K$ by a fundamental system of solutions of $\mathcal{L}(Y)=0$ and with constant field equal to $C$. In the case $C$ algebraically closed, it is known that the Picard-Vessiot extension exists and is unique up to $K$-differential isomorphism (see \cite{kol2}). In \cite{kolex}, Kolchin quotes a remark of Baer who notes that the difficulty lies not in proving the existence of a fundamental system of solutions of the given differential equation but in proving the existence of one which brings in no new constants.

In \cite{sei}, Seidenberg constructed an example of a linear differential equation defined over a differential field $K$ with constant field the field $\R$ of real numbers, for which no Picard-Vessiot extension exists (see Example \ref{seid}). At first sight, this example seems to indicate that it is not possible to obtain a general result on existence of Picard-Vessiot extension beyond the class of differential fields with algebraically closed field of constants (see \cite{bor} Section 5.3). Some misinterpretation of this example, quoted by several specialists, may explain the fact that such a general result has not been obtained since now. However, the differential field $K$ in Seidenberg's example is not a real field (see Definition \ref{real}). In this paper we present an existence theorem of Picard-Vessiot extensions for real differential fields with real closed field of constants and establish a Galois correspondence theorem for these Picard-Vessiot extensions. From this result, we expect to obtain a characterization of linear differential equations solvable by real Liouville functions, a question raised by O. Gel'fond and A. Khovanskii (see \cite{GK}, Remark 3. It is worth noting that the field of rational functions $\R(x_1,\dots,x_n)$ and the field of real meromorphic functions are real fields hence our work may also lead to applications in real analytic mechanics (see \cite{mor}, \cite{aud}).

It is well known that if $L$ is a splitting field of a polynomial in $K[X]$, the extension $L|K$ is normal, i.e. for any $a \in L\setminus K$, there exists $\sigma \in Aut_K L$ such that $\sigma(a) \neq a$. If $L|K$ is a normal algebraic extension, then, for any field $F$ with $K\subset F \subset L$, $L|F$ is normal as well. In his quest for a good concept of normality for differential field extensions, Kolchin observed that the direct analog of normality for differential field extensions is defective, as the property does not translate to intermediate differential fields (see Example \ref{weak}). He defined then a differential field extension $L|K$ to be normal when for any differential field $F$ with $K\subset F \subset L$ and any $a \in L\setminus F$, there exists a differential automorphism $\sigma$ of $L$ over $F$ such that $\sigma(a) \neq a$. However, the Galois correspondence theorem for normal differential extensions has some failures. Kolchin finally introduced the concept of strongly normal extension (see Definition \ref{nor}) and obtained a satisfactory Galois correspondence theorem for this class of extensions without assuming the field of constants of the differential base field to be algebraically closed (see \cite{kol} Chapter VI). Note that, for a strongly normal extension $L|K$, in the case when the constant field of $K$ is not algebraically closed, the differential Galois group is no longer the group $DAut_K L$ of $K$-differential automorphisms of $L$, rather one has to consider as well $K$-differential morphisms of $L$ in larger differential fields. Strongly normal extensions have been studied by several authors after Kolchin (see \cite{kov}, \cite{pil}, \cite{ume}). Note that, except Umemura, they assume the field of constants of the base differential field to be algebraically closed. In fact, Umemura introduces the more general concept of automorphic extension.

It is worth noting that a Picard-Vessiot extension of a differential field with algebraically closed field of constants is normal, in Kolchin's sense, but this is no longer true for a Picard-Vessiot extension of a real differential field with real closed field of constants. However a Picard-Vessiot extension is always strongly normal. In the case of Picard-Vessiot extensions of real fields with real closed field of constants, we can adopt a definition of the differential Galois group inspired by Kolchin's but simpler than his one. We obtain then a Galois correspondence theorem which classifies intermediate differential fields of a Picard-Vessiot extension of a real differential field with real closed field of constants in terms of its differential Galois group.

In this paper, we shall deal with ordinary differential
fields of characteristic 0. We shall denote by $C_{K}$ the constant field of the differential field $K$.

We refer the reader to \cite{boch} for topics on real field theory, to \cite{haj}, \cite{mag} or \cite{ps} for topics on differential Galois theory. We shall
use the terminology of Kolchin on constrained extensions \cite{kol}, \cite{kol3}.

We thank Ehud Hrushovski, Marius van der Put, Michael Singer and Marcus Tressl for interesting discussions along the elaboration of this paper.

\section{Preliminaries}

We recall now the notions of normality for differential field extensions introduced by Kolchin  and the precise definition of Picard-Vessiot extension. We adopt Kovacic's definition of strong morphism and strongly normal extension (see \cite{kov}).

\begin{definition} \label{nor} Let $L|K$ be an extension of differential fields.

\begin{enumerate}
\item $L|K$ is \emph{weakly normal} if for every $a \in L \setminus K$, there exists $\sigma \in DAut_K L$ such that $\sigma(a) \neq a$.
\item $L|K$ is \emph{normal} if for every differential field $F$, with $K\subset F \subset L$, $L|F$ is weakly normal.
\item If $M$ is a differential field extension of $L$, $f:L \rightarrow M$ is a differential $K$-morphism, we say that $f$ is \emph{strong} if the following two conditions are satisfied.
    \begin{enumerate}
    \item $f(a)=a$, for all $a \in C_L$,
    \item the following equalities between composite fields hold.

    $$L f(L)=LC(f)=f(L)C(f),$$

    \noindent
     where $C(f)$ is the field of constants of the composite field $L f(L)$.
\end{enumerate}
\item $L|K$ is \emph{strongly normal} if  it is finitely differentially generated and for any differential field extension $M$ of $L$, any $K$-differential morphism $f$  of $L$ in $M$ is strong.
\end{enumerate}
\end{definition}

\begin{definition}\label{PV} Given a homogeneous linear differential equation

$$\mathcal{L}(Y):=Y^{(n)}+a_{n-1}Y^{(n-1)}+ \ldots + a_{1}Y'+a_{0}Y=0$$

\noindent of
order $n$ over a differential field $K$, a differential extension
$L|K$ is a \emph{Picard-Vessiot
extension} for $\cL$ if
\begin{enumerate}
\item $L=K\langle \eta_1,\dots,\eta_n\rangle$, where $\eta_1,\dots,\eta_n$ is a fundamental
set of solutions of $\cL(Y)=0$  in $L$.
\item Every constant of $L$ lies in $K$, i.e. $C_K=C_L$.
\end{enumerate}
\end{definition}

As mentioned in the introduction, in the case when the constant field $C_K$ of the differential field $K$ is algebraically closed, it is known that there exists a Picard-Vessiot extension for a given
homogeneous linear ordinary differential
equation defined over $K$ which is unique, up to $K$-differential isomorphism. The following example due to Seidenberg (\cite{sei}) proves that one cannot expect a Picard-Vessiot extension to exist for any linear differential equation over an arbitrary differential field.

\begin{example}\label{seid}
We consider the field of real numbers $\R$ with trivial derivation and the differential field $K$ obtained by adjoining to $\R$ a solution of the differential equation
$4a^2+a'^2=-1$, such that $a' \neq 0$. Let us look at the homogeneous linear differential equation $Y''+Y=0$ defined over $K$. Seidenberg proved that for any differential field extension $L$ of $K$ containing a solution of this last equation, the inclusion of $\R$ in the constant field of $L$ is strict. In other words, there is no Picard-Vessiot extension of $K$ for this equation.
\end{example}

In this paper we shall deal with linear differential equations defined over real differential fields with real closed field of constants. We recall now the meaning of real and real closed field  and some of their properties.

\begin{definition}\label{real}
An \emph{ordered field} is a field endowed with an ordering compatible with the field operations.
A field $K$ is called a \emph{real field} if it can be ordered or equivalently if $-1$ is not a sum of squares in $K$. A real field $K$ which has no nontrivial real algebraic extensions is called a \emph{real closed field}.
An algebraic extension $L$ of an ordered field $K$ is called a \emph{real closure} of $K$ if $L$
is real closed and the inclusion $K \hookrightarrow L$ preserves the ordering of $K$.
\end{definition}

A real field always has characteristic zero. If $K$ is a real field, the ring $K[i]:=K[X] / (X^2+1)$ is a field which is a quadratic extension of $K$. If $K$ is a real field, the field of rational functions $K(X)$ is as well real.

A field $K$ is a real closed field if and only if the ring $K[i]$ is an algebraically closed field.

Every ordered field $K$ has a real closure which is unique up to $K$-isomorphism.
The fields $\Q$ and $\R$ with their natural orderings are clearly real fields. Moreover $\R$ is a real closed field.

Let us note that the field $K$ in example \ref{seid} is not real since, by construction, $-1$ is a sum of squares in $K$. However the class of real differential fields with real closed field of constants will be a good setting to establish the existence of Picard-Vessiot extensions. We note that a partial result in this direction has been obtained in \cite{sow}. Recently, H. Gillet, S. Gorchinskiy and A. Ovchinnikov \cite{GGO} proved the existence of Picard-Vessiot extensions for real differential fields with real closed field of constants using Tannakian categories. However the fact that the Picard-Vessiot extension is real  does not follow from their result.

\vspace{0.5cm}
Regarding Galois correspondence theorem, it is worth noting that a Picard-Vessiot extension of a differential field with algebraically closed field of constants is normal, in Kolchin's sense (see Definition \ref{nor}). This is not longer true for Picard-Vessiot extensions of real fields with real closed field of constants, though, as it can be seen in the following example.

\begin{example}\label{weak}
Let us consider the real differential field $K:=\R(t)$, with derivation $\frac{d}{dt}$. Its constant field is clearly $\R$. We consider the differential field extension $L:=K(e^{t})|K$. It is a Picard-Vessiot extension for the equation $Y'=Y$, defined over $K$. The assignment $e^t \mapsto \lambda e^t$ determines a $K$-differential automorphism of $L$, for all $\lambda \in \R^*$, so $L|K$ is weakly normal. Now consider the intermediate field $F=K(e^{3t})$. The only $F$-automorphism of $L$ is identity, hence $L|F$ is not weakly normal, so $L|K$ is not normal.
\end{example}

\section{Existence of Picard-Vessiot extensions for real fields}

In this section we shall establish the existence of a real Picard-Vessiot extension for a homogeneous linear differential equation defined
over a real differential field. Our result was announced in \cite{CRAS}. For the reader's convenience we give its proof here. We recall the notions of constrained extension and constrainedly closed field of Kolchin and some of their properties.

\begin{definition} An element $\eta$ in an extension of $\F$ is
said to be \emph{constrained} over $\F$ if there exists a
differential polynomial $C \in \F\{y\}$ with $C(\eta)\neq 0$ such
that $C(\eta')=0$, for every non-generic differential
specialization $\eta'$ of $\eta$ over $\F$.

We shall say that the differential field extension $\G \supset\F$ is \emph{constrained} if and only if
each element of $\G$ is constrained over $\F$.
\end{definition}

\begin{remark} Under the hypothesis that $\F$ contains nonconstant elements, our definition of constrained extension coincides with Kolchin's one (cf. \cite{kol3}, section 2).
\end{remark}

The proof of the following proposition is given in \cite{kol}, Chap. III, Section 10, Proposition 7.

\begin{proposition}\label{const} Every algebraic extension is constrained.
Moreover, if $\F$ contains nonconstant elements, the field of constants of a constrained extension of
$\F$ is algebraic over the field of constants $C_{\F}$ of $\F$.
\end{proposition}

We are interested in differential fields which are a quadratic extension of a real differential subfield.
We set our terminology in the following definition.

\begin{definition}
\begin{enumerate}
\item A pair of differential fields $(\F,K)$ is of \emph{real type} if
$K$ is a formally real field and $\F=K(i)$, where $i$ is a
root of $X^2+1$.
\item If $(\F,K)$ is a pair of differential fields of real type, a \emph{differential extension of real type of $(\F,K)$}
 is a  pair of differential fields $(\G,L)$ of real type such that $\F \subset \G$ and $K= L\cap \F$.
\end{enumerate}
\end{definition}

M. Singer has proved that for a real differential
field $K$, there exists a differential extension $(\Phi,\K)$ of real type of $(K(i),K)$ such that
$\Phi$ is a constrainedly closed field and $\K$ is a real closed field (cf. \cite{Si},
Theorem).

\begin{theorem}\label{th} Let $(\F,K)$ be a pair of differential fields of real type. There exists a differential extension of real type $(\G,L)$ of $(\F,K)$ such that $\G$ is constrainedly closed and a constrained
extension of $\F$.
\end{theorem}

\noindent {\it Proof.} Fix a differential extension $(\Phi,\K)$ of real type of $(\F,K)$ such that
$\Phi$ is a constrainedly closed field and $\K$ is a real closed field (Singer's pair) and let $\M$ be the family of differential extensions $(\F',K')$ of real type of $(\F,K)$ contained in $(\Phi,\K)$ such that
$\F'$ is a constrained extension of $\F$. Then
$\M\neq \emptyset$ and is partially ordered by inclusion.

As $\M$ is ordered inductively, by Kuratowski-Zorn's lemma $\M$ has a maximal
element. Let us fix such a maximal element $(\G,L)$. We claim
that $\G$ is constrainedly closed and $(\G,L)$ is of real type.

To establish this fact, consider any element
$\eta \in \Phi$ that is constrained over $\G$, say with
constraint $B$ and let us observe that the conjugate element $\overline{\eta}$ is constrained over $\G$ with constraint $\overline{B}$ (where we consider the conjugation $c$ in $\Phi$ determined by $c_{\K}=Id_{\K}$ and $c(i)=-i$). Therefore $(\G\langle \eta, \overline{\eta}\rangle,\G\langle \eta, \overline{\eta}\rangle\cap \K)$ is an element of $\M$. By the maximality of $(\G,L)$,
then $\G\langle \eta, \overline{\eta}\rangle=\G$. In particular,
$\eta \in \G$. Because of \cite{kol3} Section 3, Corollary 2, this
establishes our claim.

\begin{theorem}\label{teo}  Let $(\F,K)$ be a pair of differential fields of real type and let $(\G,L)$ be a differential extension of real type of $(\F,K)$ such that $\G$ is a constrainedly closed and constrained extension of $\F$. Let $\cL(Y) \in \F\{Y\}$ be a homogeneous linear differential polynomial of order $n$. There exist $y_1,\dots,y_n \in \G$ solutions of $\cL(Y)=0$ such that $y_1,\dots,y_n$ are linearly independent over $C_{\G}$.
\end{theorem}

\noindent {\it Proof.} We denote by $Wr(y_1,\dots,y_n)$ the wronskian of $y_1,\dots,y_n$. For $y_1,\dots,y_m$ with $m<n$ and $Wr(y_1,\dots,y_{m})\neq 0$, $Wr(y_1,\dots,y_{m+1})$ has order $m$ and because $\G$ is constrainedly closed we can find $y_{m+1} \in \G$ such that $\cL(y_{m+1})=0$ and $Wr(y_1,\dots,y_{m+1})\neq 0$. Therefore we can find a system of solutions for $\cL(Y)=0$ linearly independent over $C_{\G}$.

\begin{corollary}\label{cor} Let $K$ be a real differential field with real closed field of constants $C_{K}$. Let $\cL(Y)=0$ be a homogeneous linear differential equation defined over $K$. Then there exists a Picard-Vessiot extension $L$ of $K$ for the equation $\cL(Y)=0$ and moreover $L$ is a real field.
\end{corollary}

\noindent {\it Proof.} Let $\F=K(i)$ and $(\G,L)$ be the differential extension of real type of $(\F,K)$ such that $\G$ is constrainedly closed and a constrained
extension of $\F$, given by theorem \ref{th}. We assume that $\F$ contains nonconstant elements. By Proposition \ref{const}, we have $C_{\G}=C_{\F}$. Let us denote by $c$ the conjugation of $\G$ determined by $c(i)=-i, c_{|L}=Id_L$. Let $V$ be the $C_{\F}$-subspace of $\G$ generated by the $C_{\F}$-linearly independent solutions $y_1,\dots,y_n$ of $\cL(Y)=0$, given by Theorem \ref{teo}, and let $V^c$ be the $C_{\F}$-subspace of $V$ fixed by the conjugation $c$. The differential subfield of $\G$ generated by $K=\F^c$ and $V^c$ is a real Picard-Vessiot extension of $K$ for the equation $\cL(Y)=0$.

 If $K=C_K$, we may consider the field of rational functions $K(t)$ and extend derivation by $t'=1$. Then, if $\cL(Y)=0$ is a homogeneous linear differential equation defined over $K$, there exists a Picard-Vessiot extension $L$ of $K(t)$ for $\cL(Y)=0$, which is a real field. Then the subfield of $L$ differentially generated over $K$ by the $K$-vector space of solutions for $\cL(Y)=0$ in $L$ is a real Picard-Vessiot extension of $K$ for $\cL(Y)=0$.

\begin{remark} In the case of a linear differential equation defined over a differential field $\F$ with algebraically closed constant field, the Picard-Vessiot extension is proved to be unique, up to $\F$-differential isomorphism. Let us consider a pair $(\F,K)$ of differential fields of real type such that the constant field $C$ of $K$ is real closed, a linear differential equation $\mathcal{L}(Y)=0$, defined over $K$ and a Picard-Vessiot extension $L$ of $K$ for $\mathcal{L}(Y)=0$. The set of $K$-isomorphism classes of Picard-Vessiot extensions of $K$ for $\mathcal{L}(Y)=0$ is in bijection with $H^1(\Gal(\F|K),DAut_{\F}(\G))$, where $\G=L(i)$. So, the uniqueness of the Picard-Vessiot extension does not hold in the case of differential fields with a real closed field of constants. For example, $H^1(\Gal(\C|\R),SO(n,\R))$ is not trivial, as it can be identified with the set of equivalence classes of quadratic forms of rank $n$ with positive discriminant (see \cite{ser} III 3.2). However, if we moreover assume the differential field to be real and we want to restrict to real Picard-Vessiot extensions, the problem of uniqueness is more subtle (see example \ref{sin} below). It is connected to the problem of determining the isomorphism classes of real fields $K$ having isomorphic extensions $K(i)$ which is, as far as we know, not solved.

One can construct examples of linear differential equations defined over a real differential field with real closed field of constants having more than one real Picard-Vessiot extension (see example \ref{and} below). It would be interesting to determine if the uniqueness of the Picard-Vessiot extension holds for a fixed ordering of the real differential base field.
\end{remark}

\begin{example}\label{sin} Let us consider the differential equation $Y''+Y=0$ defined over the field $\F=\C(t)$, with derivation $d/dt$.  Its Picard-Vessiot extension is $\G=\F(\sin t, \cos t)$ and its differential Galois group is $\SO(2,\C)$.

We consider now the same equation over the field $K=\R(t)$. We have two Picard-Vessiot extensions of $K$ for this equation which are not $K$-isomorphic, namely $L_1=K(\sin t, \cos t)$ and $L_2=K(i\sin t, i\cos t)$ corresponding to the two elements in $H^1(\Gal(\C|\R),\SO(2,\C))$. We observe that $L_1$ is a real field, while $L_2$ is not as $(i\sin t)^2+(i\cos t)^2=-1$.
\end{example}

\begin{example}\label{and} The quadratic field extensions $\R(\sqrt{t})|\R(t)$ and $\R(\sqrt{-t})|\R(t)$ are Picard-Vessiot extensions for the equation $Y'=Y/(2t)$ over $\R(t)$. Both fields $\R(\sqrt{t})$ and $\R(\sqrt{-t})$ are real. In the first, we have $t>0$, while in the second, we have $t<0$. Now, $t \mapsto -t$ defines an isomorphism between these two fields but they are not $\R(t)$-isomorphic.
\end{example}

\section{Galois correspondence}\label{g}

As mentioned in the introduction, a Picard-Vessiot extension is strongly normal. Hence, the fundamental theorem established by Kolchin in \cite{kol} chapter VI applies to Picard-Vessiot extensions. However, for a strongly normal extension $L|K$, Kolchin defines the differential Galois group by means of differential $K$-isomorphisms of $L$ in the differential universal extension of $L$. Moreover, he uses Weil language of algebraic geometry. In this section, we give a more direct definition of the differential Galois group of a Picard-Vessiot extension over a differential field with real closed field of constants, we endow it with a linear algebraic group structure and establish a Galois correspondence theorem in our setting.

\subsection{Galois group}

Let $K$ be a real differential field with real closed field of constants $C$, $\F=K(i)$. For a  real Picard-Vessiot extension $L|K$, we shall consider the set $DHom_K(L,\G)$ of $K$-differential morphisms from $L$ into $\G=L(i)$. We shall see that we can define a group structure on this set and we shall take it as the differential Galois group $\DGal(L|K)$ of the Picard-Vessiot extension $L|K$. We shall prove that it is a $C$-defined (Zariski) closed subgroup of some $\overline{C}$-linear algebraic group, where $\overline{C}$ denotes the algebraic closure of $C$.

We observe that we can define mutually inverse bijections

$$\begin{array}{ccc} DHom_K(L,\G) & \rightarrow & DAut_{\F}\G \\ \sigma & \mapsto & \widehat{\sigma}  \end{array}, \quad \begin{array}{ccc}  DAut_{\F}\G & \rightarrow & DHom_K(L,\G) \\ \tau & \mapsto & \tau_{|L} \end{array},$$

\noindent where $\widehat{\sigma}$ is the extension of $\sigma$ to $\G$ defined by $\widehat{\sigma}(a+i\,b)=\sigma(a)+i \sigma(b)$, for $a,b \in L$. We may then transfer the group structure from
$DAut_{\F}\G$ to $DHom_K(L,\G)$.

Let now $\eta_1,\dots,\eta_n$ be $C$-linearly independent elements in $L$ such that \linebreak $L=K\langle \eta_1,\dots,\eta_n\rangle$ and $\sigma \in DHom_K(L,\G)$. We have then $\sigma(\eta_j)=\sum_{i=1}^n c_{ij} \eta_i,$ \linebreak $1\leq j \leq n$, with $c_{ij} \in \overline{C}$. We may then associate to $\sigma$ the matrix $(c_{ij})$ in $\GL(n,\overline{C})$.

\begin{proposition} Let $K$ be a real differential field with real closed field of
constants $C$, $L=K\langle \eta_1,\dots,\eta_n\rangle$ a real Picard-Vessiot
extension of $K$, where $ \eta_1,\dots,\eta_n$ are $C$-linearly independent, $\F=K(i), \G=L(i)$. There exists a set $S$ of polynomials
$P(X_{ij}), 1\leq i,j \leq n$, with coefficients in $C$ such that
\begin{enumerate}[1)]
\item If $\sigma \in DHom_K(L,\G)$ and
$\sigma(\eta_j)=\sum_{i=1}^n c_{ij}\eta_i$, then $P(c_{ij})=0, \forall P
\in S$.
\item Given a matrix $(c_{ij}) \in \GL(n,\overline{C})$ with $P(c_{ij})=0, \forall P
\in S$, there exists a differential $K$-morphism $\sigma$ from
$L$ to $\G$ such that $\sigma(\eta_j)=\sum_{i=1}^n c_{ij}\eta_i$.
\end{enumerate}
\end{proposition}

\noindent {\it Proof.} The proof follows the steps of prop. 6.2.1 in \cite{haj}.
Let $K\{Z_1,\dots,Z_n\}$ be the ring of
differential polynomials in $n$ indeterminates over $K$. We define
a differential $K$-morphism $\varphi$ from $K\{Z_1,\dots,Z_n\}$ in $L$ by
$Z_j \mapsto \eta_j$. Then $\Gamma:=\Ker \varphi$ is a prime differential
ideal of $K\{Z_1,\dots,Z_n\}$. Let $\G[X_{ij}], 1\leq i,j \leq n$,
be the ring of polynomials in the indeterminates $X_{ij}$ with the
derivation defined by $X_{ij}'=0$. We define a differential
$K$-morphism $\psi$ from $K\{Z_1,\dots,Z_n\}$ to $\G[X_{ij}]$ such that
$Z_j\mapsto \sum_{i=1}^n X_{ij} \eta_i$. Let $\Delta:=\psi(\Gamma)$. Let $\{w_k\}$ be a basis of the
$C$-vector space $\G$. We write each polynomial in $\Delta$ as a
linear combination of the $w_k$ with coefficients polynomials in
$C[X_{ij}]$. We take $S$ to be the collection of all these
coefficients.

\noindent 1. Let $\sigma$ be a differential $K$-morphism from
$L$ to $\G$ and $\sigma(\eta_j)=\sum_{i=1}^n c_{ij} \eta_i$. We consider the
diagram

\begin{picture}(250,150)(0,100)
\put(110,200){$K\{Z_1,\dots,Z_n\}$}
\put(190,204){\vector(1,0){50}}
\put(210,210){$\varphi$}
\put(250,200){$L$}
\put(150,190){\vector(0,-1){50}}
\put(155,165){$\psi$}
\put(253,190){\vector(0,-1){50}}
\put(140,120){$\G[X_{ij}]$}
\put(250,120){$\G$}
\put(260,165){$\sigma$}
\put(180,123){\vector(1,0){60}}
\put(210,130){$v$}
\put(87,220){$Z_j$}
\put(110,219){\line(0,1){8}}
\put(110,223){\vector(1,0){130}}
\put(250,220){$\eta_j$}
\put(88,210){\line(1,0){8}}
\put(92,210){\vector(0,-1){75}}
\put(70,120){$\sum X_{ij}\eta_i$}
\put(142,100){$X_{ij}$}
\put(180,102){\vector(1,0){60}}
\put(180,106){\line(0,-1){8}}
\put(247,100){$c_{ij}$}
\end{picture}

\vspace{0.4cm}
 \noindent It is clearly commutative. The image of
$\Gamma$ by $\sigma \circ \varphi$ is
$0$. Its image by $v\circ \psi$ is $\Delta$ evaluated in $X_{ij}=c_{ij}$. Therefore
all polynomials of $\Delta$ vanish at $c_{ij}$. Writing this down
in the basis $\{w_k\}$, we see that all polynomials of $S$ vanish
at $c_{ij}$.

\noindent 2. Let us now be given a matrix $(c_{ij}) \in \GL(n,\overline{C})$
such that $F(c_{ij})=0$ for every $F$ in $S$. We consider the differential morphism

$$ \begin{array}{ccc} K\{Z_1,\dots,Z_n\}&\rightarrow &
\G \\ Z_j & \mapsto & \sum_i c_{ij} \eta_i
\end{array}.$$

\noindent By the
hypothesis on $(c_{ij})$, and the definition of the set $S$, we
see that the kernel of this morphism contains $\Gamma$ and so, we
have a differential $K$-morphism

$$ \begin{array}{cccc} \sigma: & K\{\eta_1,\dots,\eta_n\}&\rightarrow &
\G \\ & \eta_j & \mapsto & \sum_i c_{ij} \eta_i
\end{array}.$$

\noindent Taking into account that the elements $\sigma(\eta_j), 1\leq j \leq r$, are $\overline{C}$-linearly independent, we obtain that $\sigma$ is injective and so extends to a $K$-differential morphism from $L$ to $\G$.
\hfill $\Box$

\vspace{0.5cm}
If $L|K$ is a real Picard-Vessiot extension for a homogeneous linear differential equation of order $n$ defined over $K$, the preceding proposition gives that $\DGal(L|K)$ is a $C$-defined closed subgroup of $\GL(n,\overline{C})$.

\vspace{0.5cm}
Real Picard-Vessiot extensions satisfy the following normality property.

\begin{proposition}
Let $K$ be a real differential field with real closed field of constants~$C$, $L|K$ a real differential Picard-Vessiot extension. For $a \in L\setminus K$, there exists a $K$-differential morphism $\sigma:L \rightarrow \G$ such that $\sigma(a) \neq a$.
\end{proposition}

\noindent {\it Proof.} As $\G|\F$ is a Picard-Vessiot extension and the constant field $\overline{C}$ of $\F$ is algebraically closed, we know (\cite{haj} prop. 6.1.2) that there exists an $\F$-differential automorphism $\widehat{\sigma}$ of $\G$ such that $\widehat{\sigma}(a) \neq a$. We can then take $\sigma=\widehat{\sigma}_{|L}$. \hfill $\Box$

\vspace{0.5cm}
For a subset $S$ of $\DGal(L|K)$, we set $L^S:=\{ a \in L : \sigma(a)=a, \, \forall \sigma \in S \}$.

\begin{corollary}\label{cor}
Let $K$ be a real differential field with real closed field of constants $C$, $L|K$ a real differential Picard-Vessiot extension. We have \linebreak $L^{\DGal(L|K)}=K$.
\end{corollary}

\subsection{Fundamental theorem}\label{tf}

Let $K$ be a real differential field with real closed field of constants $C$ and $L|K$ a real Picard-Vessiot extension. For a closed subgroup $H$ of $\DGal(L|K)$, $L^H$ is a differential subfield of $L$ containing $K$. If $E$ is an intermediate differential field, i.e. $K\subset E \subset L$, then $L|E$ is a real Picard-Vessiot extension and $\DGal(L|E)$ is a $C$-defined closed subgroup of $\DGal(L|K)$.

\begin{theorem} Let $L|K$ be a real Picard-Vessiot extension, $\DGal(L|K)$ its differential Galois
group.
\begin{enumerate}[1.]
\item
The correspondences

$$H\mapsto L^H \quad , \quad E\mapsto \DGal(L|E)$$

\noindent define inclusion inverting mutually inverse bijective
maps between the set of $C$-defined closed subgroups $H$ of $\DGal(L|K)$
and the set of differential fields $E$ with $K\subset E \subset
L$.

\item The intermediate differential field $E$ is a Picard-Vessiot extension of $K$ if and only if the subgroup $\DGal(L|E)$ is normal in $\DGal(L|K)$. In this case, the restriction morphism

    $$\begin{array}{ccc} \DGal(L|K)& \rightarrow & \DGal(E|K) \\ \sigma & \mapsto & \sigma_{|E} \end{array} $$

\noindent induces an isomorphism $$\DGal(L|K)/\DGal(L|E) \simeq \DGal(E|K).$$
\end{enumerate}
\end{theorem}

\noindent {\it Proof.} 1. It is clear that both maps invert inclusion. If $E$ is an intermediate differential field of $L|K$, we have $L^{\DGal(L|E)}=E$, taking into account that $L|E$ is Picard-Vessiot and corollary \ref{cor}. For $H$ a $C$-defined closed subgroup of $\DGal(L|K)$, the equality $H=\DGal(L|L^{H})$ follows from the correspondent equality in Picard-Vessiot theory for differential fields with algebraically closed field of constants (\cite{haj} theorem 6.3.8).

\noindent 2. If $E$ is a Picard-Vessiot extension of $K$, then  $E(i)$ is a Picard-Vessiot extension of $K(i)$ and so $\DGal(L|E)$ is normal in $\DGal(L|K)$. Reciprocally, if $\DGal(L|E)$ is normal in $\DGal(L|K)$, then the subfield of $L(i)$ fixed by $\DGal(L|E)$ is a Picard-Vessiot extension of $K(i)$. Now, this field is $E(i)$. So, if $E(i)$ is differentially generated over $K(i)$ by a $\overline{C}$-vector space of finite dimension $V$, then $E$ is differentially generated over $K$ by the $C$-vector space $V^{c}=\{y \in V: c(y)=y\}$, where $c$ is the conjugation of $E(i)$ determined by $c(i)=-i$, $c_{|E}=Id_E$. Hence $F|K$ is a real Picard-Vessiot extension. Let us note that $V$ may be assumed to be $c$-stable since if $E(i)=K(i)\langle V \rangle$, then also $E(i)=K(i)\langle V+\overline{V}\rangle$, where $\overline{V}=\{c(v): v \in V\}$. The last statement of the theorem follows from the fundamental theorem of Picard-Vessiot theory in the case of algebraically closed fields of constants (\cite{haj} theorem 6.3.8). \hfill $\Box$

\begin{remark} All results in section \ref{g} remain valid for $K$ any differential field with real closed field of constants $C$ and $L|K$ a Picard-Vessiot extension. Just observe that, as $-1$ is not a square in $C$, $K(i)$ is also in this case a quadratic extension of $K$. However, without assuming $K$ real, we cannot assure that a Picard-Vessiot extension exists for a given linear differential equation defined over $K$.
\end{remark}

\vspace{1cm}
\footnotesize
\begin{tabular}{lcl}
Teresa Crespo && Zbigniew Hajto \\
Departament d'\`{A}lgebra i Geometria && Faculty of Mathematics and Computer Science \\  Universitat de Barcelona && Jagiellonian University \\ Gran Via de les Corts Catalanes 585 && ul. Prof. S. \L ojasiewicza 6 \\
08007 Barcelona, Spain && 30-348 Krak\'ow, Poland \\
teresa.crespo@ub.edu && zbigniew.hajto@uj.edu.pl \\
&& \\

El\.zbieta Sowa && \\
Faculty of Applied Mathematics && \\
AGH University of Science and Technology && \\
al. Mickiewicza 30 && \\
30-059 Krak\'ow,  Poland  && \\
esowa@agh.edu.pl && 
\end{tabular}

\end{document}